\documentclass{article}
\usepackage{graphicx}
\usepackage{amsmath}
\usepackage{amssymb}

\newtheorem{theorem}{Theorem}
\newtheorem{remark}{Remark}
\newtheorem{proposition}{Proposition}

\begin{document}

\title{Nash equilibria strategies and equivalent single-objective optimization problems. The case of linear partial differential equations}

\author{Angel Manuel Ramos \\
             Interdisciplinary Mathematics Institute \& \\
             Department of Applied Mathematics and Mathematical Analysis \\
             Complutense University of Madrid \\
              Spain\\
              angel@mat.ucm.es
}

\maketitle

\begin{abstract}
In this paper we study the existence and uniqueness of Nash equilibria (solution to {\em competition-wise} problems, with several controls trying to reach possibly different goals) associated to linear partial differential equations and show that, in some cases, they are also the solution of suitable single-objective optimization problems (i.e. {\em cooperative-wise} problems, where all the controls cooperate to reach a common goal). We use cost functions associated with a particular linear parabolic partial differential equations and distributed controls, but the results are also valid for more general linear differential equations (including elliptic and hyperbolic cases) and controls (e.g. boundary controls, initial value controls,...).
\end{abstract}

\vspace{.2cm}


\section{Introduction}

Nash equilibria are solutions of a noncooperative multiobjective optimization
strategy first proposed by Nash (see \cite{Nash}). Since it originated in game
theory and economics, the notion of player is often used. For an optimization
problem with $N$ objectives or functionals $J_i$ to minimize, a Nash strategy
consists in having $N$ players or controls $v_i$, each optimizing his own
criterion. However, each player has to optimize his criterion given that all
the other criteria are fixed by the rest of the players. When no player can
further improve his criterion, it means that the system has reached a Nash
equilibrium state.

Of course, there are other strategies for multiobjective optimization,
such as the Pareto cooperative strategy (\cite{Pareto}), the Stackelberg hierarchical strategy (\cite{Stack}) or the Stackelberg-Nash strategy (\cite{DiLi}).

To the best of our knowledge, \cite{RaGlPe1} and \cite{RaGlPe2} are the first articles dealing with the theoretical and numerical study of Nash equilibria for differential games associated to partial differential equations. Following \cite{RaGlPe1} we deal here with a general linear case with $N$ cost functions and controllers and show how, in some cases, the Nash equilibria (solution to differential games associated to multiobjective optimization problems with several noncooperative controllers), are also the solution of single-objective optimization problems (where all the controllers cooperate to reach a common goal). We use cost functions associated with linear parabolic partial differential equations and distributed controls, but other kinds of linear differential equations (e.g. elliptic, hyperbolic,..) and controls (e.g. boundary controls, initial value controls,...) can be also used, using the same technique.

The fact that a noncooperative game (i.e a competition-wise problem) can be seen as a (cooperative) single-objective optimization problem (i.e. a noncompetition-wise problem) is very interesting, not only because of the curious noncooperative-cooperative equivalence, but also because of the huge amount of software to compute solutions and literature written about the latter kind of problems, that could be used in the framework of, apparently, a different type of problems.

In Section \ref{fotp} we formulate the problem and give an optimality system providing a necessary and sufficient condition for the Nash equilibria. The existence and uniqueness of Nash equlibria is studied in Section \ref{eauof}.  In Section \ref{espcas} we show the equivalence, in some cases, between the noncooperative multiobjective differential games defining the Nash equlibria and suitable (cooperative) single-objective optimization problems. Finally, in Section \ref{conclus} give a summary of the major results of the paper.

\section{Formulation of the Problem} \label{fotp}

Let us consider $T>0$, $\Omega\subset \mathbb{R}^d$ a bounded and smooth open set with $d\in \{ 1,2,3\}$, and two
subsets $\Gamma_1$, $\Gamma_2\subset \partial \Omega$ such that
$\partial\Omega =\Gamma_1 \cup \Gamma_2$. We define $Q=\Omega\times
(0,T)$,
$\Sigma_1=\Gamma_1\times (0,T)$, $\Sigma_2=\Gamma_2\times (0,T)$ and
the control Hilbert spaces  ${\cal U}_i=L^2(\omega_i \times (0,T))$ and ${\mathcal U}={\cal U}_1\times \cdots \times {\cal U}_N$, where $N\in \mathbb{N}$, $\omega_i\subset\Omega$ , for $i\in \{ 1,...,N\}$, and $\omega_i\cap\omega_j
=\emptyset$ if $i\neq j$. Finally, we consider the functionals $J_i:{\mathcal U}\rightarrow \mathbb{R}$, with $i\in \{ 1,...,N\}$, given by
\begin{eqnarray*}
J_i(v_1,...,v_N) & = & \frac{\alpha_i}{2}\int_{\omega_i\times (0,T)} |v_i|^2{\rm d}x {\rm d}t
\\
& & +\frac{1}{2}\int_{Q} \rho_i(x)
|y-y_{i,{\rm d}}|^2 {\rm d}x{\rm d}t+ \frac{1}{2}\int_{\Omega}
\eta_i(x)|y(T)-y_{i,{\rm T}}|^2 {\rm d}x,
\end{eqnarray*}
for every $v=(v_1,...,v_N)\in {\cal U}$, where  $\alpha_i >0$, $\rho_i , \eta_i \in L^\infty (\Omega)$ such that $\rho_i,
\eta_i\geq 0$, the function $y=y(v)$ is defined as the
solution of
\begin{equation} \label{eqdy}
\left\{
\begin{array}{ll}
{\displaystyle \frac{\partial y}{\partial t} - \Delta y=f+\sum_{i=1}^N v_i\chi_{\omega_i}} & \mbox{ in  }Q, \\
y(0)=y_0 & \mbox{ in }\Omega , \\
{\displaystyle y=g_1} & \mbox{ on }\Sigma_1 , \\
{\displaystyle \frac{\partial y}{\partial n}=g_2} & \mbox{ on }\Sigma_2 , \\
\end{array}
\right.
\end{equation}
with $f,g_1,g_2,y_0,y_{i,{\rm d}}$ and $y_{i,{\rm T}}$ being smooth enough functions and $\chi_\omega:\Omega\rightarrow \mathbb{R}$ the characteristic function (with values 1 in $\omega$ and 0 in $\Omega\setminus \omega$) for any $\omega\subset\Omega$.

This generalizes the typical examples in the literature of 2 controls (instead of $N$), $\rho_i =k_i\chi\omega_{{\rm d},i}$ and $\eta_i =l_i \chi\omega_{{\rm T},i}$, where $k_i,l_i>0$, $\omega_{{\rm d},i}, \omega_{{\rm T},i}\subset \Omega$. A special case is when $\omega_{{\rm T}1}\cap \omega_{{\rm T},2}\neq \emptyset$ and\verb+/+or
$\omega_{{\rm d},1}\cap \omega_{{\rm d}2} \neq \emptyset$. This case is a {\em competition-wise} problem, with each control (or {\em player}) trying to reach (possibly) different goals over a common {\em domain}. In some sense this is the case where the behavior of the solution $y$ associated to a Nash equilibrium is most difficult to forecast.

\begin{remark}
Most of the results to follow are also valid for more general linear operators
such as, for instance,
$${\cal A}\varphi = \frac{\partial \varphi}{\partial t}- \nabla \cdot
(A(x)\nabla
\varphi) +V\cdot \nabla \varphi +c(x)\varphi .$$
The technique is also valid for different type of controls such as, for instance,
boundary or initial controls. $\square$
\end{remark}

Now, given $i\in \{ 1,...,N\}$, for every $(w_1,...,w_{i-1},w_{i+1},...,w_N)\in {\cal U}_1\times \cdots \times {\cal U}_{i-1}\times {\cal U}_{i+1} \times \cdots \times {\cal U}_N$ we consider the optimal control
problem
$$({\cal CP}_i)\left\{
\begin{array}{l}
\mbox{Find }
u_i(w_1,...,w_{i-1},w_{i+1},...,w_N)\in {\cal U}_i, \mbox{ such that} \\
J_i(w_1,...,w_{i-1},u_i(w_1,...,w_{i-1},w_{i+1},...,w_N),w_{i+1},...,w_N) \\
\hspace*{3cm} \leq J_1(w_1,...,w_{i-1},v_i,w_{i+1},...w_N), \ \forall v_i\in {\cal U}_i.
\end{array}
\right.
$$

The (unique) solution $u_i(w_1,...,w_{i-1},w_{i+1},...,w_N)$ of problem $({\cal
CP}_i)$ is characterized by
$$
\frac{\partial J_i}{\partial v_i}(w_1,...,w_{i-1}, u_i(w_1,...,w_{i-1},w_{i+1},...,w_N), w_{i+1},...,w_N)=0.
$$

\noindent
Therefore, a Nash equilibrium is a $N$-tuple $(u_1,...,u_N)\in {\cal U}$ such
that $u_i=u_i(u_1,...u_{i-1},u_{i+1},...,u_N)$ for all $i\in \{ 1,...,N\}$, i.e. $(u_1,...,u_N)$ is a solution of the
{\em
coupled (optimality) system}:
\begin{equation} \label{j1j2}
\frac{\partial J_i}{\partial v_i} (u_1,...,u_N)=0, \ \ \ \ \forall \ i\in\{ 1,...,N\}.
\end{equation}
In the linear case studied here, this system of equations is a necessary and sufficient condition for $u$ to be a Nash equilibrium. In general this system is only a necessary condition, although in some nonlinear cases (see, e.g. \cite{RaRo}), the functionals are convex and system (\ref{j1j2}) is also a sufficient condition.

Following \cite{RaGlPe1} it is easy to prove that, if $i\in \{ 1,...,N\}$,
$$
\frac{\partial J_i}{\partial v_i}(v)=\alpha_i v_i+p_i(v)\chi_{\omega_i} \in {\cal U}_i,
$$
where for any $v=(v_1,...,v_N)\in {\cal U}$ the function $p_i=p_i(v)$ is the solution of the adjoint system
$$
\left\{
\begin{array}{ll}
{\displaystyle -\frac{\partial p_i}{\partial t} -\Delta p_i
=\rho_i(y-y_{i,{\rm d}})} &
\mbox{ in }Q, \\
p_i(T)=\eta_i(y(T)-y_{i,{\rm T}}) & \mbox{ in }\Omega , \\
p_i=0 & \mbox{ on }\Sigma_1 , \\
{\displaystyle \frac{\partial p_i}{\partial n}=0} & \mbox{ on }\Sigma_2
\end{array}
\right.
$$
and $y=y(v)$ is the solution of (\ref{eqdy}).

Therefore, system (\ref{j1j2}) is equivalent to the (optimality) system
$$
\left\{
\begin{array}{l}
{\displaystyle u_i=-\frac{1}{\alpha_i}p_i\chi_{\omega_i}, \ \ i\in \{1,...,N\}} \\
\
\\
\left\{
\begin{array}{ll}
{\displaystyle \frac{\partial y}{\partial t} - \Delta y=f+\sum_{i=1}^N u_i\chi_{\omega_i} } & \mbox{ in  }Q, \\
y(0)=y_0 & \mbox{ in }\Omega , \\
{\displaystyle y=g_1} & \mbox{ on }\Sigma_1 , \\
{\displaystyle \frac{\partial y}{\partial n}=g_2} & \mbox{ on }\Sigma_2 ; \\
\end{array}
\right.
\\
\ \\
\left\{
\begin{array}{ll}
{\displaystyle -\frac{\partial p_i}{\partial t} -\Delta p_i
=\rho_i(y-y_{i,{\rm d}})} &
\mbox{ in }Q, \\
p_i(T)=\eta_i(y(T)-y_{i,{\rm T}}) & \mbox{ in }\Omega , \\
p_i=0 & \mbox{ on }\Sigma_1 , \\
{\displaystyle \frac{\partial p_i}{\partial n}=0} & \mbox{ on }\Sigma_2  .
\end{array} \nonumber
\right. , \ \ i\in \{1,...,N\}
\end{array}
\right.
$$

\section{Existence and uniqueness of solution of Nash Equlibria} \label{eauof}

It is obvious that
\begin{equation} \label{am}
 v=(v_1,...,v_N)
\longrightarrow
(\frac{\partial J_1}{\partial v_1}(v_1,...,v_N) ,...,
\frac{\partial J_N}{\partial v_N}(v_1,...,v_N))
\in {\cal U}
\end{equation}
is an affine mapping of ${\cal U}$. Therefore,
there exist a linear continuous mapping ${\cal A}\in {\cal L}({\cal U},{\cal U})$
and a vector $b\in {\cal U}$ such that
$$(\frac{\partial J_1}{\partial v_1}(v_1,...,v_N) ,...,
\frac{\partial J_N}{\partial v_N}(v_1,...,v_N))
={\cal A} v -b.
$$
Let us identify mapping ${\cal A}$: For every $v=(v_1,...,v_N)\in {\cal U}$, the
linear part of the affine mapping in relation (\ref{am}) is defined by
$${\cal A}v= (\alpha_1 v_1 +\widetilde{p}_1\chi_{\omega_1} ,..., \alpha_N v_N
+\widetilde{p}_N\chi_{\omega_N}),
$$
where $\widetilde{p}_i=\widetilde{p}_i(v)$, $i \in\{ 1,...,N\}$, is the solution of
$$
\left\{
\begin{array}{ll}
{\displaystyle -\frac{\partial \widetilde{p}_i}{\partial t} -\Delta \widetilde{p}_i
=\rho_i\widetilde{y}} &
\mbox{ in }Q, \\
\widetilde{p}_i(x,T)=\eta_i\widetilde{y}(T) & \mbox{ in }\Omega , \\
\widetilde{p}_i=0 & \mbox{ on } \Sigma_1 , \\
{\displaystyle \frac{\partial \widetilde{p}_i}{\partial n}=0} & \mbox{ on }\Sigma_2 ,
\end{array}
\right.
$$
and $\widetilde{y}=\widetilde{y}(v)$ is the solution of
$$
\left\{
\begin{array}{ll}
{\displaystyle \frac{\partial \widetilde{y}}{\partial t} - \Delta \widetilde{y}=\sum_{i=1}^Nv_i\chi_{\omega_i}} & \mbox{ in  }Q, \\
\widetilde{y}(0)=0 & \mbox{ in }\Omega , \\
{\displaystyle \widetilde{y}=0} & \mbox{ on }\Sigma_1 , \\
{\displaystyle \frac{\partial \widetilde{y}}{\partial n}=0} & \mbox{ on }\Sigma_2 . \\
\end{array}
\right.
$$
\begin{proposition} \label{sp}
Mapping ${\cal A}:{\cal U}\rightarrow {\cal U}$ is linear and continuous. Furthermore, if ${\displaystyle \min_{i\in\{ 1,...,N\}} \{ \alpha_i\}}$ is sufficiently large, it is also ${\cal U}$-elliptic, i.e., there existe $C>0$ such that
$$ ({\cal A}v , v)\geq C ||v||^2,$$
where $(\cdot,\cdot)$ and $||\cdot||$ represent the canonical scalar product and norm of the Hilbert space ${\cal U}$, respectively.
\end{proposition}
{\bf Proof:} It is obvious that ${\cal A}$ is a linear mapping and it is easy
to show that it is continuous (see \cite{LiMa}).

Let us consider $v=(v_1,...,v_N)\in {\cal U}$ and $w=(w_1,...,w_N)\in {\cal U}$. We have then
$$({\cal A}v , w)
= \big( (\alpha_1 v_1 +\widetilde{p}_1\chi_{\omega_1} ,..., \alpha_N v_N
+\widetilde{p}_N\chi_{\omega_N}) , (w_1,...,w_N) \big)
$$
$$
=\sum_{i=1}^N\int_{\omega_i\times (0,T)} (\alpha_i v_i +\widetilde{p}_i(v))w_i{\rm d}x{\rm d}t.
$$
Let us focus on the term $\int_{\omega_i\times (0,T)}
\widetilde{p}_i(v)w_i{\rm d}x{\rm d}t$, following the approach in \cite{RaGlPe1} and \cite{CaFe19}. We have
$$\int_{\omega_i\times (0,T)} \widetilde{p}_i(v)w_i{\rm d}x{\rm d}t
$$
$$=\int_{Q}
\widetilde{p}_i(v) \left(\frac{\partial}{\partial t} \widetilde{y}(0,...,w_i,...,0)-\Delta
\widetilde{y}(0,...,w_i,...,0)\right){\rm d}x{\rm d}t
$$
$$=\int_{Q}\left( -\frac{\partial}{\partial t} \widetilde{p}_i(v)-\Delta
\widetilde{p}_i(v)\right) \widetilde{y}(0,...,w_i,...,0){\rm d}x{\rm d}t
$$
$$+\int_{\Omega} \eta_i\widetilde{y}(T;v)\widetilde{y}(T;0,...,w_i,...,0){\rm d}x $$
$$=\int_{Q}\rho_i\widetilde{y}(v)\widetilde{y}(0,...,w_i,...,0){\rm d}x{\rm d}t
+\int_{\Omega} \eta_i \widetilde{y}(T;v)\widetilde{y}(T;0,...,w_i,...,0){\rm d}x.$$
Then,
$$
({\cal A}v,w)
 =  \sum_{i=1}^N \left( \alpha_i \int_{\omega_i\times (0,T)} v_iw_i{\rm d}x{\rm d}t \right.
 $$
$$
\left. + \int_{Q}\rho_i \widetilde{y}(v)\widetilde{y}(0,...,w_i,...,0){\rm d}x{\rm d}t
+\int_{\Omega} \eta_i\widetilde{y}(T;v)\widetilde{y}(T;0,...,w_i,...,0){\rm d}x\right) . $$
Since the mapping $v\rightarrow \widetilde{y}(v)$ is linear and continuous from ${\cal U}$ to ${\mathcal C}([0,T];L^2(\Omega))$ (see, e.g., \cite{LiMa}), it is easy to prove there exist a constant $c>0$ such that
$||\widetilde{y}(v)||_{L^2(Q)}+||\widetilde{y}(T;v)||_{L^2(\Omega)}\leq c || v||$.
Therefore,
\begin{eqnarray*}
({\cal A}(v,v) & \geq & \min_{i\in \{ 1,...,N\}} \{ \alpha_i\}||v||^2 - \sum_{i=1}^N c^2(||\rho_i||_{L^\infty(\Omega)}+||\eta_i||_{L^\infty(\Omega)}) ||v|| \ ||v_i||_{{\cal U}_i} \\
& \geq & C||v||^2 ,
\end{eqnarray*}
with
$$
C=\left(\min_{i\in \{ 1,...,N\}} \{ \alpha_i\}-\sum_{i=1}^N c^2(||\rho_i||_{L^\infty(\Omega)}+||\eta_i||_{L^\infty(\Omega)}) \right) ||v||^2.
$$
Notice that $C>0$ if ${\displaystyle \min_{i\in\{ 1,...,N\}} \{ \alpha_i\}>\sum_{i=1}^N c^2(||\rho_i||_{L^\infty(\Omega)}+||\eta_i||_{L^\infty(\Omega)})}$, which proves that ${\cal A}$ is ${\cal U}$-elliptic in that case and completes the
proof. $\blacksquare$

\vspace{.2cm}

Let us identify $b$: The constant part of the affine mapping (\ref{am}) is
the function $b\in {\cal U}$ defined by $b=(\overline{p}_1\chi_{\omega_1},...,\overline{p}_N\chi_{\omega_N}),$
where $\overline{p}_i$, $i\in \{1,...,N\}$, is the solution of

$$\left\{
\begin{array}{ll}
{\displaystyle -\frac{\partial \overline{p}_i}{\partial t} -\Delta \overline{p}_i
=\rho_i(\overline{y}-y_{i,{\rm d}})} &
\mbox{ in }Q, \\
\overline{p}_i(T)=\eta_i(\overline{y}-y_{i,{\rm T}}) & \mbox{ in }\Omega , \\
\overline{p}_i=0 & \mbox{ on }\Sigma_1 , \\
{\displaystyle \frac{\partial \overline{p}_i}{\partial n}=0} & \mbox{ on }\Sigma_2 ,
\end{array}
\right.
$$
and $\overline{y}$ is the solution of
\begin{equation} \label{edyb}
\left\{
\begin{array}{ll}
{\displaystyle \frac{\partial \overline{y}}{\partial t} - \Delta \overline{y}=f} & \mbox{ in  }Q, \\
\overline{y}(0)=y_0 & \mbox{ in }\Omega , \\
{\displaystyle \overline{y}=g_1} & \mbox{ on }\Sigma_1 , \\
{\displaystyle \frac{\partial \overline{y}}{\partial n}=g_2} & \mbox{ on }\Sigma_2 . \\
\end{array}
\right.
\end{equation}
Notice that, for any $v\in {\cal U}$, $y(v)=\widetilde{y}(v)+\overline{y}$ and $p_i(v)=\widetilde{p}_i(v)+\overline{p}_i$.

\begin{theorem} \label{eunw1} If ${\displaystyle \min_{i\in\{ 1,...,N\}} \{ \alpha_i\}}$ is sufficiently large, there exist a unique Nash equilibrium of the problem defined in Section \ref{fotp}.
\end{theorem}
{\bf Proof:} As showed above, the Nash equilibria are characterized by the solutions of (\ref{j1j2}), which are also characterized by the solutions $u\in {\cal U}$ of
$$a(u,v)=L(v), \ \ \forall \ v\in {\cal U},
$$
where a $a(\cdot ,\cdot):{\cal U}\times {\cal U}\rightarrow \mathbb{R}$ is defined by
$$a(v,w)=({\cal A}v,w)\ \ \ \forall \ v,w\in {\cal U},$$
and $L:{\cal U}\rightarrow \mathbb{R}$ by
$$L(v)=(b,v), \ \ \ \forall \ v\in {\cal U}.$$
Proposition \ref{sp} proves that mapping $a(\cdot ,\cdot)$ is bilinear,
continuous and, if ${\displaystyle \min_{i\in\{ 1,...,N\}} \{ \alpha_i\}}$ is sufficiently large, it is also ${\cal U}$-elliptic. Furthermore, mapping $L$ is (obviously) linear and
continuous. Thus, by the (well-known) Lax-Milgram Theorem, system (\ref{j1j2}) has a unique solution or, equivalently, there exists a unique Nash equilibrium of the problem defined in Section \ref{fotp}, if ${\displaystyle \min_{i\in\{ 1,...,N\}} \{ \alpha_i\}}$ is sufficiently large. $\blacksquare$

The discretization of the problem considered above and the development of suitable algorithms to get a numerical solution approximating the Nash equilibra can follow the approaches in \cite{RaGlPe1} and \cite{CaFe19}.

\section{Equivalent single-objective control problems} \label{espcas}

In this section we will show that, in some cases, the solution of noncooperative differential games defining Nash equilibria, are the solution of suitable optimization problems, where all the controls cooperate to minimize a suitable single-objective cost function.

Let us consider the subfamily of problems defined in Section \ref{fotp}, for which $\rho_i=\rho$ and $\eta_{i}=\eta$, for all $i\in \{ 1,...,N\}$. Therefore, in this case
the functional $J_i$, with $i\in \{ 1,...,N\}$, is given by
\begin{eqnarray*}
J_i(v_1,...,v_N) & = & \frac{\alpha_i}{2}\int_{\omega_i\times (0,T)} |v_i|^2{\rm d}x {\rm d}t
\\
& & +\frac{1}{2}\int_{Q}
\rho |y-y_{i,{\rm d}}|^2 {\rm dx}{\rm d}t + \frac{1}{2}\int_{\Omega}
\eta |y(T)-y_{i,{\rm T}}|^2 {\rm d}x,
\end{eqnarray*}
with $y=y(v)$ being the solution of (\ref{eqdy}).

As in the general case studied in Section \ref{fotp}, a Nash equilibrium is a $N$-tuple $(u_1,...,u_N)\in {\cal U}={\cal U}_1\times \cdots \times {\cal U}_N$ solution of (\ref{j1j2}), where
$$
\frac{\partial J_i}{\partial v_i}(v)=\alpha_i v_i+p_i(v)\chi_{\omega_i} \in {\cal U}_i, \ \ \forall i\in\{ 1,...,N\}
$$
for any $v=(v_1,...,v_N)\in {\cal U}$ and $p_i=p_i(v)$ is now the solution of
$$
\left\{
\begin{array}{ll}
{\displaystyle -\frac{\partial p_i}{\partial t} -\Delta p_i
=\rho (y-y_{i,{\rm d}})} &
\mbox{ in }Q, \\
p_i(T)=\eta (y(x,T)-y_{i,{\rm T}}) & \mbox{ in }\Omega , \\
p_i=0 & \mbox{ on }\Sigma_1 , \\
{\displaystyle \frac{\partial p_i}{\partial n}=0} & \mbox{ on }\Sigma_2.
\end{array}
\right.
$$

Therefore, system (\ref{j1j2}) is equivalent in this case to
$$
\left\{
\begin{array}{l}
{\displaystyle u_i=-\frac{1}{\alpha_i}p_i\chi_{\omega_i}, \ \ i\in \{1,...,N\}} \\
\
\\
\left\{
\begin{array}{ll}
{\displaystyle \frac{\partial y}{\partial t} - \Delta y=f+\sum_{i=1}^N u_i\chi_{\omega_i} } & \mbox{ in  }Q, \\
y(0)=y_0 & \mbox{ in }\Omega , \\
{\displaystyle y=g_1} & \mbox{ on }\Sigma_1 , \\
{\displaystyle \frac{\partial y}{\partial n}=g_2} & \mbox{ on }\Sigma_2 ; \\
\end{array}
\right.
\\
\ \\
\left\{
\begin{array}{ll}
{\displaystyle -\frac{\partial p_i}{\partial t} -\Delta p_i
=\rho(y-y_{i,{\rm d}})} &
\mbox{ in }Q, \\
p_i(T)=\eta (y(T)-y_{i,{\rm T}}) & \mbox{ in }\Omega , \\
p_i=0 & \mbox{ on }\Sigma_1 , \\
{\displaystyle \frac{\partial p_i}{\partial n}=0} & \mbox{ on }\Sigma_2  .
\end{array} \nonumber
\right. , \ \ i\in \{1,...,N\}
\end{array}
\right.
$$
Again
$$(\frac{\partial J_1}{\partial v_1}(v_1,...,v_N) ,...,
\frac{\partial J_N}{\partial v_N}(v_1,...,v_N))
={\cal A} v -b.
$$
and now
$${\cal A}v= (\alpha_1 v_1 +\widetilde{p}\chi_{\omega_1} ,..., \alpha_N v_N
+\widetilde{p}\chi_{\omega_N}),
$$
where $\widetilde{p}=\widetilde{p}(v)$ is the solution of
$$
\left\{
\begin{array}{ll}
{\displaystyle -\frac{\partial \widetilde{p}}{\partial t} -\Delta \widetilde{p}
=\rho \widetilde{y}} &
\mbox{ in }Q, \\
\widetilde{p}(x,T)=\eta \widetilde{y}(T) & \mbox{ in }\Omega , \\
\widetilde{p}=0 & \mbox{ on } \Sigma_1 , \\
{\displaystyle \frac{\partial \widetilde{p}}{\partial n}=0} & \mbox{ on }\Sigma_2 ,
\end{array}
\right.
$$
\begin{proposition} \label{sp2} For the family of problems studied in Section \ref{espcas}
mapping ${\cal A}:{\cal U}\rightarrow {\cal U}$ is linear, continuous self-adjoint and ${\cal U}$-elliptic.
\end{proposition}
{\bf Proof:} Following the proof of Proposition \ref{sp}, ${\cal A}$ is a linear and continous mapping. Furthermore, given $v=(v_1,...,v_N)\in {\cal U}$ and $w=(w_1,...,w_N)\in {\cal U}$,
$$
({\cal A}v,w)
 =  \sum_{i=1}^N \Big( \alpha_i \int_{\omega_i\times (0,T)} v_iw_i{\rm d}x{\rm d}t
 $$
$$
+ \int_{Q}\rho \widetilde{y}(v)\widetilde{y}(0,...,w_i,...,0){\rm d}x{\rm d}t
+\int_{\Omega} \eta \widetilde{y}(T;v)\widetilde{y}(T;0,...,w_i,...,0){\rm d}x\Big)
$$
$$
=  \sum_{i=1}^N \alpha_i \int_{\omega_i\times (0,T)} v_iw_i{\rm d}x{\rm d}t
+ \int_{\Omega}\rho \widetilde{y}(v)\widetilde{y}(w){\rm d}x{\rm d}t
+\int_{\Omega} \eta \widetilde{y}(T;v)\widetilde{y}(T;w){\rm d}x
$$
$$
=(v,{\cal A}w)
$$
and
$$
({\cal A}v,v)\geq \min_{i\in\{ 1,...,N\}} \{ \alpha_i\} ||v||^2,
$$
which proves that ${\cal A}$ is self-adjoint and ${\cal U}$-elliptic. $\blacksquare$

\vspace{.2cm}

The constant part of the affine mapping (\ref{am}) is
the function $b\in {\cal U}$ defined by $b=-(\overline{p}_1\chi_{\omega_1},...,\overline{p}_N\chi_{\omega_N}),$
where $\overline{p}_i$, $i\in \{1,...,N\}$, is now the solution of

$$\left\{
\begin{array}{ll}
{\displaystyle -\frac{\partial \overline{p}_i}{\partial t} -\Delta \overline{p}_i
=\rho (\overline{y}-y_{i,{\rm d}})} &
\mbox{ in }Q, \\
\overline{p}_i(T)=\eta (\overline{y}-y_{i,{\rm T}}) & \mbox{ in }\Omega , \\
\overline{p}_i=0 & \mbox{ on }\Sigma_1 , \\
{\displaystyle \frac{\partial \overline{p}_i}{\partial n}=0} & \mbox{ on }\Sigma_2 ,
\end{array}
\right.
$$
and $\overline{y}$ is the solution of (\ref{edyb}).

\begin{theorem} There exist a unique Nash equilibrium of the problem defined in Section \ref{espcas}.
\end{theorem}
{\bf Proof:} The proof follows the one of Theorem \ref{eunw1}, taking into account that in this case ${\cal A}$ is unconditionally ${\cal U}$-elliptic. $\blacksquare$

The discretization of the problem considered above and the development of suitable algorithms to get a numerical solution approximating the Nash equilibra are given in \cite{RaGlPe1}, where numerical examples are also showed.
\begin{theorem} The (unique) Nash equilibrium $u=(u_1,...,u_N)\in {\cal U}$ of the problem defined in Section \ref{espcas} is the (unique) solution of the following optimal control problems:
\begin{itemize}
\item Find $u\in {\cal U}$ such that ${\displaystyle J(u)=\min_{v\in V}J(v)}$, where
\begin{eqnarray*}
J(v) & = & \sum_{i=1}^N \frac{\alpha_i}{2}\int_{\omega_i\times (0,T)} |v_i|^2{\rm d}x {\rm d}t \\
& & +\frac{1}{2}\sum_{i=1}^N\left( \int_{Q} \rho
|y(0,...,v_i,...,0)-y_{i,{\rm d}}|^2 {\rm dx}{\rm d}t \right. \\
& & \left.\hspace*{3cm} + \int_{\Omega} \eta |y(T; 0,...,v_i,...,0)-y_{i,{\rm T}}|^2 {\rm d}x \right) \\
& &  +2\sum_{i,j=1 (i<j)}^N\left( \int_{Q} \rho y(0,...,v_i,...,0) y(0,...,v_j,...,0) {\rm dx}{\rm d}t \right.\\
& &  \left.\hspace{2cm} \int_{\omega_{T}} \eta y(T;0,...,v_i,...,0)y(T;0,...,v_j,...,0){\rm d}x\right) .
\end{eqnarray*}
\item Given $j,p\in \{ 1,...,N\}$, find $u\in {\cal U}$ such that ${\displaystyle J_{j,p}(u)=\min_{v\in {\cal U}}J_{j,p}(v)}$, where
\begin{eqnarray*}
J_{j,p}(v) & = & \sum_{i=1}^N \frac{\alpha_i}{2}\int_{\omega_i\times (0,T)} |v_i|^2{\rm d}x {\rm d}t \\
& & +\frac{1}{2}\int_{Q} \rho
|y(v)-y_{j,{\rm d}}|^2 {\rm dx}{\rm d}t + \frac{1}{2}\int_{\Omega} \eta |y(T;v)-y_{p,{\rm T}}|^2 {\rm d}x \\
& &  +2\sum_{i=1 (i\neq j)}^N\int_{Q} \rho (y_{j,{\rm d}}-y_{i,{\rm d}})\widetilde{y}0,...,v_j,...,0) {\rm dx}{\rm d}t \\
& & +2\sum_{i=1 (i\neq p)}^N \int_{\Omega} \eta (y_{p,{\rm T}}-y_{i,{\rm T}})\widetilde{y} (T;0,...,v_j,...,0) {\rm dx}.
\end{eqnarray*}
\end{itemize}
\end{theorem}
{\bf Proof:} We have seen previously that there exists a unique Nash equilibrium, which is the solution $u\in {\cal U}$ of
$$
({\cal A}u,v)= (b,v) \ \ \forall \ v\in {\cal U}.
$$
Then, because of the properties of ${\cal A}$ given in Theorem \ref{sp2}, we have (see, e.g., \cite[Theorem 2.44]{RaFEM}) that
$$
\widetilde{J}(u)=\min_{v\in {\cal U}} \widetilde{J}(v),
$$
with
$$
\widetilde{J}(v) = ({\cal A}v,v)-2 (b,v)= \sum_{i=1}^N \alpha_i \int_{\omega_i\times (0,T)} |v_i|^2{\rm d}x{\rm d}t
+ \int_{Q}\rho |\widetilde{y}(v)|^2{\rm d}x{\rm d}t
$$
$$
+\int_{\Omega} \eta |\widetilde{y}(T;v)|^2{\rm d}x+2 \sum_{i=1}^N \int_{Q} \rho (\overline{y}-y_{i,{\rm d}})\widetilde{y}(0,...,v_i,...,0){\rm d}x{\rm d}t
$$
$$
+2\sum_{i=1}^N \int_{\Omega} \eta (\overline{y}(T)-y_{i,{\rm T}})\widetilde{y}(T;0,...,v_i,...,0){\rm d}x.
$$
Then, using that $\widetilde{y}(v)={\displaystyle \sum_{i=1}^N \widetilde{y}(0,...,v_i,...,0)}$, we have that
$$
\widetilde{J}(v)= \sum_{i=1}^N \left( \alpha_i \int_{\omega_i\times (0,T)} |v_i|^2{\rm d}x{\rm d}t \right.
$$
$$
+ \int_{Q}\rho \left[ |\widetilde{y}(0,...,v_i,...,0)|^2 +2 (\overline{y}-y_{i,{\rm d}})\widetilde{y}(0,...,v_i,...,0) \right] {\rm d}x{\rm d}t
$$
$$
\left. \int_{\Omega} \eta \left[ |\widetilde{y}(T;0,...,v_i,...,0)|^2 +2 (\overline{y}(T)-y_{i,{\rm T}})\widetilde{y}(T;0,...,v_i,...,0)\right] {\rm d}x \right)
$$
$$
+2\sum_{i,j=1 (i<j)}^N\left( \int_{Q} \rho \widetilde{y}(0,...,v_i,...,0) \widetilde{y}(0,...,v_j,...,0) {\rm dx}{\rm d}t\right.
$$
$$
\left. \hspace{2cm} +\int_{\Omega} \eta \widetilde{y}(T;0,...,v_i,...,0) \widetilde{y}(T;0,...,v_j,...,0) {\rm dx} \right) .
$$
Hence, using that $y=\widetilde{y}+\overline{y}$ we have that
$$
\widetilde{J}(v)=J(v)-C,
$$
where
$$
C= \sum_{i=1}^N \left( \rho\int_{Q} (\overline{y}-y_{i,{\rm d}})^2 {\rm d}x{\rm d}t +\int_{\Omega} \eta (\overline{y}-y_{i,{\rm T}})^2 {\rm d}x \right)
$$
is a constant (independent of $v$), which completes the proof of the first part of the theorem.

\vspace{.2cm}

In order to prove the second part of the theorem, given $j\in \{ 1,...,N\}$, let us focus on the following terms of $\widetilde{J}(v)$:
$$
\int_{Q}\rho |\widetilde{y}(v)|^2{\rm d}x{\rm d}t +2 \sum_{i=1}^N \int_{Q} \rho (\overline{y}-y_{i,{\rm d}})\widetilde{y}(0,...,v_i,...,0){\rm d}x{\rm d}t
$$
$$
= \int_{Q}\rho \left[ |\widetilde{y}(v)|^2 +2 (\overline{y}-y_{j,{\rm d}})( \widetilde{y}(v)-\widetilde{y}(v_1,...,v_{j-1},0,v_{j+1},..., v_N)\right] {\rm d}x{\rm d}t
$$
$$
+2 \sum_{i=1 (i\neq j)}^N \int_{Q} \rho (\overline{y}-y_{i,{\rm d}})\widetilde{y}(0,...,v_i,...,0){\rm d}x{\rm d}t
$$
$$
= \int_{Q}\rho \left[ |\widetilde{y}(v)|^2 +2 (\overline{y}-y_{j,{\rm d}})\widetilde{y}(v)\right] {\rm d}x{\rm d}t
$$
$$
+2 \sum_{i=1 (i\neq j)}^N \int_{Q} \rho (y_{j,{\rm d}}-y_{i,{\rm d}})\widetilde{y}(0,...,v_i,...,0){\rm d}x{\rm d}t.
$$
Something similar can be done with other terms of $\widetilde{J}(v)$, so that
$$
\widetilde{J}(v)=J(v)-C_{j,p},
$$
where
$$
C_{j,p}=  \int_{Q} \rho (\overline{y}-y_{j,{\rm d}})^2 {\rm d}x{\rm d}t +\int_{\Omega} \eta (\overline{y}-y_{p,{\rm T}})^2 {\rm d}x ,
$$
which completes the proof. $\blacksquare$

\section{Conclusions} \label{conclus}

This paper studies Nash equilibria of noncooperative differential games with several players (controllers), each one trying to minimize his own cost function defined in terms of a general class of linear partial differential equations. We give results of existence and uniqueness of Nash equilibria and show how, in some cases, the corresponding Nash equilibria (solution to {\em competition-wise} problems, with each control trying to reach possibly different goals), are also the solution of suitable single-objective optimization problems (i.e. {\em cooperative-wise} problems, where all the controls cooperate to reach a common goal). A natural question arises: Are there Nash equilibria associated to nonlinear problems than can be also characterized as the solutions of single-objective problems? This is an open problem for interested researchers.

\vspace{.3cm}
{\bf Acknowledgements:}
The research of the author was partially supported by the Spanish Ministry of Economy and Competitiveness under project MTM2015-64865-P (MINECO / FEDER), and the Research Group MOMAT (Ref. 910480) of the Complutense University of Madrid.

\end{document}